\documentclass[letterpaper, 10 pt, conference]{ieeeconf}  

\IEEEoverridecommandlockouts                              

\usepackage{graphics} 
\usepackage{epsfig} 
\usepackage{mathptmx} 
\usepackage{times} 
\usepackage{amsmath} 
\usepackage{amssymb}  
\usepackage[hidelinks]{hyperref}
\usepackage{subfig}
\usepackage[x11names]{xcolor}
\usepackage[maxbibnames=10,giveninits=true,sorting=none]{biblatex}
\newtheorem{assumption}{Assumption}
\newtheorem{proposition}{Proposition}
\newtheorem{remark}{Remark}

\title{\LARGE \bf
A Sequential Quadratic Programming Perspective on Optimal Control
}

\author{Abhijeet and Suman Chakravorty
\thanks{The authors are with the Department of Aerospace Engineering, Texas A\&M University, College Station, TX 77843 USA. \{\tt abhinir, schakrav\}@tamu.edu
}
}

\begin{document}

\maketitle
\thispagestyle{empty}
\pagestyle{empty}

\begin{abstract}
This paper offers a unified perspective on different approaches to the solution of optimal control problems through the lens of constrained sequential quadratic programming. In particular, it allows us to find the relationships between Newton's method, the iterative LQR (iLQR), and Differential Dynamic Programming (DDP) approaches to solve the problem.
It is shown that the iLQR is a principled SQP approach, rather than simply an approximation of DDP by neglecting the Hessian terms, to solve optimal control problems that can be guaranteed to always produce a cost-descent direction and converge to an optimum; while Newton's approach or DDP do not have similar guarantees, especially far from an optimum. Our empirical evaluations on the pendulum and cart-pole swing-up tasks serve to corroborate the SQP-based analysis proposed in this paper.
\end{abstract}

\section{Introduction}

The study of iterative trajectory optimization has a long history in nonlinear control. Among the most influential methods are Differential Dynamic Programming (DDP), introduced in the 1960s~\cite{mayne1965ddp}, and iterative Linear Quadratic Regulation (iLQR), which later emerged as a simplified yet powerful alternative. Early research established DDP as a Newton-type algorithm for trajectory optimization with locally quadratic convergence~\cite{mayne1965ddp, JacobsonMayne1970DDP}, whereas iLQR was designed around consistent first-order approximations that make it scalable and numerically reliable~\cite{todorov2005ilqg, tassa2012synthesis}. Over the decades, these two approaches have become foundational tools for nonlinear control, trajectory optimization, and reinforcement learning.\\

A large body of literature has analyzed their theoretical properties and limitations. Several works have emphasized the quadratic convergence of DDP under ideal assumptions~\cite{mayne1965ddp, JacobsonMayne1970DDP, Pantoja_SN}, but have also noted the difficulties in practice when the algorithm’s backward pass combines second-order expansions with first-order terms. Such inconsistencies can compromise stability and necessitate heavy regularization~\cite{Pantoja_SN}. An explicit link between DDP and Newton’s method was established in~\cite{Pantoja_SN}, and it is argued in \cite{LiaoShoemaker1992} that, for discrete-time optimal control, DDP can outperform Newton’s method—retaining strong local convergence while offering greater robustness and efficiency under practical problem structures—thereby clarifying both its strengths and its limits away from the optimum. Other studies further examined mixed algorithms that blend Newton-type steps with more robust strategies, claiming improved stability in practice~\cite{LiaoShoemaker1992}. For rigid-body systems, accelerated second-order DDP has been explored \cite{NgangaWensing2021RAL}. Meanwhile, extensions such as iterative Linear Quadratic Gaussian (iLQG)~\cite{todorov2005ilqg} and control-limited iLQR~\cite{tassa2014control} have reinforced the importance of iLQR in stochastic and constrained settings. Beyond classical control, iLQR-based approaches have also been integrated into policy search~\cite{levine2013guided}, legged locomotion~\cite{zhang2025wholebodymodelpredictivecontrollegged}, and optimization-based planning and feedback design~\cite{manchester2017dircol}. Related constrained/unscented formulations have also been studied \cite{Plancher2017CUDP}.\\

In contrast to these developments, the adoption of DDP has been more domain-specific. In aerospace engineering, for instance, DDP has been applied to trajectory optimization in astrodynamics and spacecraft guidance, where its second-order structure can be leveraged in high-precision planning~\cite{Russell1, Russell2, Aziz2019HDDP_CR3BP}. Beyond astrodynamics, DDP has been used for complex robotic systems~\cite{Oshin2022PDDP_RSS}, and planetary landing problems~\cite{Noyes2021RobustEntryGuidance}, demonstrating its utility in solving nonlinear, constrained, and multi-phase optimization problems. While effective in these contexts, these studies also highlight the algorithm’s sensitivity to problem structure and its dependence on strong regularization, particularly when the local quadratic models become inconsistent.\\

Albeit owing its inception to interpretation as an approximation of DDP by neglecting the second-order dynamics terms, we show that iLQR is a principled algorithm which can be derived independently of DDP and that often outperforms DDP in practice. The current work reveals a striking contrast between the two methods. DDP is theoretically appealing because of its use of second-order expansions and the possibility of quadratic convergence near an optimum owing to its similarity to Newton's method. However, the inconsistency of its backward pass can yield unreliable cost predictions, non-descent directions, or even divergence. The iLQR, by contrast, employs a consistent first-order approximation for both linear and quadratic terms, ensuring guaranteed descent directions. Although often described as a simplification of DDP, iLQR is a complete algorithm in its own right, with strong theoretical guarantees and practical advantages that frequently outweigh the appeal of DDP especially when far from an optimum.\\ 

\textit{Contributions.} In this paper, we revisit the solution of optimal control problems from both theoretical and empirical perspectives. We analyze optimal control methods within a sequential quadratic programming (SQP) framework, which provides a unifying lens for understanding algorithms for solving these problems. Our analysis is inspired, in part, by the broader nonlinear programming literature, where the role of second-order information and consistency in SQP-type algorithms has been carefully studied~\cite{boggs1995sqp}. The critical insight is that the Newton step in a constrained optimization problem is equivalent to solving a QP, which reduces to an LQR problem for optimal control. The SQP formulation allows us to derive theoretical convergence guarantees for iLQR and show that it is a principled approach to reduce the cost while satisfying the dynamical constraints. As a corollary, we also develop the exact Newton step to solve the optimal control problem, termed Newton LQR, which is different from the approach in the reference \cite{Pantoja_SN}. It is shown that the Newton LQR cannot be guaranteed to converge beacuse of the lack of a guaranteed descent direction, far from an optimum, owing to the Hessian terms. We also show that the DDP can be considered an approximation of the Newton LQR close to an optimum, and thus, has no guarantees of convergence far from an optimum. Through experiments on pendulum and cart-pole swing-up tasks, we demonstrate how these structural differences manifest in practice. The theoretical and empirical results suggest that iLQR often delivers more reliable and effective solutions than DDP while requiring far less computation owing to the absence of the Hessians of the dynamics in the computations, thus showing that in the case of solving optimal control problems, ``less is more". In contrast to our prior work \cite{Wang2025SearchFeedbackRL} that studied the convergence of iLQR, this work proposes a unified SQP based framework to study algorithms to solve optimal control problems, and allows us to carefully understand the relationships between the Newton's method, DDP and iLQR for solving optimal control problems, theoretically as well as empirically. Most importantly, it provides critical insight into which of these should be preferred and why?\\

The remainder of the paper is organized as follows. Section~II introduces the preliminaries and problem formulation. Section~III develops an optimal control perspective using Sequential Quadratic Programming (SQP), providing the theoretical framework for our analysis. Section~IV presents numerical comparisons on benchmark problems. Section~V concludes with insights for the development of future trajectory optimization methods.

\section{Preliminaries}
Consider the discrete-time, finite-horizon optimal control problem
\begin{equation}
    \min_{\{u_t\}_{t=0}^{T-1}}
    J \;=\; \sum_{t=0}^{T-1} c(x_t,u_t) \;+\; C_{T}(x_T),
    \label{eq:ocp}
\end{equation}
subject to the nonlinear dynamics
\begin{equation}
    x_{t+1} \;=\; f(x_t,u_t), \qquad t=0,1,\dots,T-1,
    \label{eq:dynamics}
\end{equation}
with given $x_0\in\mathbb{R}^n$, states $x_t\in\mathbb{R}^n$, and controls $u_t\in\mathbb{R}^m$. The stage cost is $c:\mathbb{R}^n\times\mathbb{R}^m\!\to\!\mathbb{R}$ and the terminal cost is $C_{T}:\mathbb{R}^n\!\to\!\mathbb{R}$. We assume $f,c, C_{T}$ are twice continuously differentiable and that the control curvature is uniformly positive definite (e.g., an $R_k\succ 0$ term in $c$) so that local control updates are well-posed.

Given a nominal trajectory $\{(\overline{x}_t,\overline{u}_t)\}$ satisfying \eqref{eq:dynamics}, we denote perturbations by $\delta x_t := x_t-\overline{x}_t$ and $\delta u_t := u_t-\overline{u}_t$. Linearizations use the Jacobians
$
f_x(\overline{x}_t,\overline{u}_t), ~
f_u(\overline{x}_t,\overline{u}_t),
$
and cost derivatives are written compactly as
$c_{x_t},c_{u_t},c_{x_t,x_t},c_{x_t,u_t},c_{u_t,u_t}$ and
$C_{T,x},C_{T,xx}$, all evaluated along the nominal.

We emphasize the structural distinction that drives the subsequent analysis: iLQR employs a first-order model of the dynamics together with a quadratic model of the cost, whereas DDP mixes second-order state information with first-order control terms. Casting the local step induced by \eqref{eq:ocp}–\eqref{eq:dynamics} as a quadratic program with linearized dynamics constraints provides a unified SQP viewpoint. In the next section, we use this QP/SQP formulation to make precise how these modeling choices govern descent directions, step acceptance, and convergence guarantees.

\section{Optimal Control: an SQP Perspective}
In the following, we study the solution to an optimal control problem from an SQP perspective utilizing the unique features of the optimal control problem. However, first, we develop some basic SQP background for the application to optimal control.

\subsection{SQP Foundations}
Consider the following constrained optimization problem:
\[
\min_{x} \; c(x) 
\quad \text{s.t.} \quad h(x) = 0.
\]

The necessary conditions are obtained by forming the augmented Lagrangian:
\[
\mathcal{L} = c + \lambda h
\]
and setting 
\[
\frac{\partial \mathcal{L}}{\partial x} = 0 
\quad \text{and} \quad 
\frac{\partial \mathcal{L}}{\partial \lambda} = 0,
\]
yielding:
\[
c_x + \lambda^* h_x = 0,\;
h(x^*) = 0,
\]
at a minimum (stationary point) \((x^*, \lambda^*)\).

In order to solve this problem, one can use first- or second-order methods. In the following, we consider the second-order Newton method to solve the problem.

Suppose that \((\bar{x}, \bar{\lambda})\) is the current solution to the problem. We then expand the Lagrangian 
\(\mathcal{L}\) to second order around the nominal \((\bar{x}, \bar{\lambda})\) to yield:
\begin{align*}
\mathcal{L}(x, \lambda) = c + \lambda h 
= \big( \bar{c} + \bar{\lambda}\,\bar{h} \big) 
+ (c_x \delta x + \bar{\lambda} h_x \delta x) \\
+ \frac{1}{2} \delta x'\big( c_{xx} + \bar{\lambda} h_{xx} \big) \delta x 
+  \delta x' \, h_x \delta \lambda 
+ \bar{h}' \, \delta \lambda\\
\equiv \bar{\mathcal{L}} + \Delta \mathcal{L},
\end{align*}
where $\bar{c} = c(\bar{x})$, $\bar{h} = h(\bar{x})$, $\bar{\mathcal{L}} = \bar{c}+\bar{\lambda}\bar{h}$ and $g_x$, $g_{xx}$ denote the gradient and Hessian of a function at the given nominal point $\bar{x}$.

In Newton's method, one optimizes the \((\Delta \mathcal{L})\) part in terms of the perturbations \((\delta x, \delta \lambda)\) and updates the solution as:
\[
x = \bar{x} + \delta x, 
\quad 
\lambda = \bar{\lambda} + \delta \lambda.
\]

The necessary conditions for obtaining the optimum of \(\Delta \mathcal{L}\) in terms of \(\delta x, \delta \lambda\) are obtained by setting
\[
\frac{\partial (\Delta \mathcal{L})}{\partial (\delta x)} = 0 
\quad \text{and} \quad 
\frac{\partial (\Delta \mathcal{L})}{\partial (\delta \lambda)} = 0,
\]
yielding:
\[
\begin{bmatrix}
c_{xx} + \bar{\lambda} h_{xx} & h_x \\
h_x & 0
\end{bmatrix}
\begin{bmatrix}
\delta x \\[6pt]
\delta \lambda
\end{bmatrix}
=
\begin{bmatrix}
- c_x - \bar{\lambda} h_x \\[6pt]
-\bar{h}
\end{bmatrix}
\]

However, taking the \(\bar{\lambda} h_x\) term to the LHS yields:
\[
\begin{bmatrix}
c_{xx} + \bar{\lambda} h_{xx} & h_x \\
h_x & 0
\end{bmatrix}
\begin{bmatrix}
\delta x \\[6pt]
\bar{\lambda} + \delta \lambda
\end{bmatrix}
=
\begin{bmatrix}
- c_x \\[6pt]
-\bar{h}
\end{bmatrix}
\]
Here, note that \(\lambda = \bar{\lambda} + \delta \lambda\) is simply the new \(\lambda\).

However, the necessary conditions above are simply the solution to the QP problem:
\begin{align}
\min_{\delta x} \; c_x \delta x + \tfrac{1}{2} \delta x' \big(c_{xx} + \bar{\lambda} h_{xx}\big) \delta x \nonumber \\
\text{s.t.} \quad h_x \delta x + \bar{h} = 0. \label{N-QP}
\end{align}

Hence, the Newton step can be accomplished by solving the above QP for the solution \((\delta x^q, \lambda^q)\) and updating the solution as:
\[
x = \bar{x} + \delta x^q, 
\quad 
\lambda = \lambda^q.
\]

Thus, the original constrained optimization can be solved by solving a sequence of QP problems above, and this will be termed as the Newton SQP approach. \\
However, in general, it is not advisable to use the ``true" Hessian of the Lagrangian \(\big(c_{xx} + \bar{\lambda} h_{xx}\big)\) in the QP unless one is close to an optimum.

To see this, note that in the constrained case, the descent direction 
$\delta X = \begin{pmatrix} \delta x \\ \delta \lambda \end{pmatrix}$ 
is given by:
\[
\delta X = (\nabla^2 \mathcal{L})^{-1} \nabla \mathcal{L},
\]
and
\[
\nabla^2 \mathcal{L} = 
\begin{bmatrix}
c_{xx} + \bar{\lambda} h_{xx} & h_x \\
h_x & 0
\end{bmatrix}
\quad \nsucc 0
\]
if the current solution is far from an optimum.  

If the Hessian $\nabla^2 \mathcal{L} \succ 0$, the Newton direction can be assured to be a descent direction, since then:
\[
\delta X' \nabla \mathcal{L} = \nabla \mathcal{L}' (\nabla^2 \mathcal{L})^{-1} \nabla \mathcal{L} > 0.
\]
However, if $\nabla^2 \mathcal{L} \nsucc 0$, there is no such guarantee.  

Thus, in order to be able to use the SQP approach, it becomes critical that one be able to ensure that the QP solution is a descent direction for the cost function.  

Consider now a modified problem (QP):
\begin{align}
\min_{\delta x} \; c_x' \delta x + \tfrac{1}{2} \, \delta x' C_{xx} \delta x \nonumber\\
h_x \delta x + \bar{h} = 0. \label{M-QP}
\end{align}

We make the assumption that $h(\bar{x}) = \bar{h} = 0$, i.e., our current estimate is feasible.  

The necessary conditions for this QP are:
\begin{align}
c_x + C_{xx} \delta x + h_x \lambda &= 0,  \\
\bar{h} + h_x \delta x &= 0. 
\end{align}

Consider the solution $\delta x$ to the necessary conditions above, then:
\[
\delta x' (c_x + C_{xx} \delta x + h_x \lambda) = 0.
\]

This implies
\[
\delta x' c_x = - \left( \delta x' C_{xx} \delta x + \delta x' h_x \lambda \right).
\]

But since $\delta x$ satisfies the constraint $\delta x' h_x = -\bar{h} = 0$, we get
\[
\delta x' c_x = - \delta x' C_{xx} \delta x < 0 
\quad \text{if } C_{xx} > 0.
\]
Since $c(x)$ is a cost function by design, it is a reasonable assumption that $C_{xx} > 0$ for all points $x$.  

Hence, the solution to the modified QP \eqref{M-QP} is always a descent direction for the cost function.  
Furthermore,
\[
h(\bar{x} + \delta x) \approx h(\bar{x}) + h_x \delta x = 0
\]
if $\delta x$ is small, and thus the update $x = \bar{x} + \delta x$ also satisfies the constraint. If $\delta x$ is not small enough, we can do a line search using a parameter $\alpha \in (0,1]$ such that
\[
h(x) = h(\bar{x} + \alpha \delta x) \approx h(\bar{x}) + \alpha h_x \delta x = 0,
\]
and thus $x = \bar{x} + \alpha \delta x$ is feasible.  

Hence, if we have a feasible solution, the modified QP always gives a descent direction for the cost while ensuring that the solution remains feasible if allied with a suitable line search approach.  
The development above can be summarized as the following result:  

\medskip

\begin{proposition} 
Consider the modified QP \eqref{M-QP}. Given $C_{xx} \succ 0$, the solution to the QP, $\delta x$, is always a descent direction for the cost function $c(x)$. Furthermore, the new solution $x = \bar{x} + \alpha \delta x$ can be made feasible, i.e., $h(x) = 0$, by a suitable choice of the line-search parameter $\alpha \in (0,1]$.
\end{proposition}

\subsection{Application to Optimal Control}

Consider the optimal control problem:
\[
\min_{\{u_t\}} \; \sum_{t=0}^{T-1} c(x_t, u_t) + c_T(x_T),
\]
subject to
\[
x_{t+1} = f(x_t, u_t).
\]

For simplicity, we shall assume the following structure for the cost and the dynamics.

\medskip

\begin{assumption}
The cost 
\[
c(x,u) = l(x) + \tfrac{1}{2} u' R u,
\]
with $R > 0$, and $l_{xx} \geq 0$ given any $x$.
\end{assumption}

\medskip

\begin{assumption}
    The dynamics are affine in the control input:
\[
x_{t+1} = \bar{f}(x_t) + \bar{g}(x_t) u_t.
\]
\end{assumption}
We do the following development for scalar state and control pair $(x,u)$ to simplify the notation, but the results can be easily generalized to the vector case.  
A final note here is that by following the dynamics, a feasible trajectory $(\bar{x}_t, \bar{u}_t)$ can always be generated.

Next, let $(\bar{x}_t, \bar{u}_t, \bar{\lambda}_t)$ be the current solution to the problem.  

Consider the following LQR problems linearized around the current $(\bar{x}_t, \bar{u}_t)$:
\begin{align}
\min_{\delta u_t}\sum_{t=0}^{T-1} 
\Big( l_{x_t} \, \delta x_t + \bar{u}_t' R \delta u_t + \, \, + \tfrac{1}{2}\delta x_t'l_{x_t,x_t}\delta x_t 
+ \tfrac{1}{2} \, \delta u_t' R \, \delta u_t \nonumber\\
+ \tfrac{\bar{\lambda}_{t+1}}{2}
\begin{bmatrix} \delta x_t \\ \delta u_t \end{bmatrix}'
\begin{bmatrix} f_{x_tx_t} & f_{x_tu_t} \\ f_{u_tx_t} & 0\end{bmatrix}
\begin{bmatrix} \delta x_t \\ \delta u_t \end{bmatrix} \Big)
+ c_{T,x} \delta x_T + \tfrac{1}{2} \, \delta x_T' c_{T,xx} \delta x_T, \nonumber\\
\delta x_{t+1} = f_{x_t} \, \delta x_t + f_{u_t} \, \delta u_t. \label{N-LQR}
\end{align}
and
\begin{align}
\min_{\delta u_t} \; \sum_{t=0}^{T-1} 
\Big( l_{x_t} \, \delta x_t + \bar{u}_t' R \, \delta u_t 
+ \tfrac{1}{2} \delta x_t' l_{x_t,x_t} \delta x_t + \,\tfrac{1}{2} \, \delta u_t' R \, \delta u_t \Big) \nonumber\\
+ c_{T,x} \delta x_T + \tfrac{1}{2} \, \delta x_T' c_{T,xx} \delta x_T,\nonumber\\
\delta x_{t+1} = f_{x_t} \, \delta x_t + f_{u_t} \, \delta u_t. \label{ILQR}
\end{align}

\medskip

Note now that the first Newton LQR problem \eqref{N-LQR} corresponds to the Newton SQP problem \eqref{N-QP}, and the second iLQR case \eqref{ILQR} corresponds to the modified SQP \eqref{M-QP}. 

The solution to the LQR problems is given by (lack of space precludes a detailed derivation but is relatively straightforward since these are simply LQR problems):
\begin{align*}
\delta x_{t+1} &= f_{x_t} \delta x_t + f_{u_t} \delta u_t, \quad \delta x_0 = 0, \\
\lambda_t &= -v_t - V_t \delta x_t, \quad \text{with } \lambda_T = -C_{T,x} - C_{T,xx} \delta x_T, \\
\delta u_t &= -k_t - K_t \delta x_t,
\end{align*}
where the solution to the Newton LQR \eqref{N-LQR} is given by:
\begin{align*}
v_t &= \ell_{x_t} + f_{x_t}' {v}_{t+1} \nonumber\\ 
- & (f_{x_t}' {V}_{t+1} f_{u_t} + \bar{\lambda}_{t+1}f_{x_tu_t})' (R+f_{u_t}' {V}_{t+1} f_{u_t})^{-1} (f_{u_t}' v_{t+1} + {R}\bar{u}_t),\\
V_t &= \ell_{x_t,x_t} + \bar{\lambda}_{t+1} f_{x_t,x_t} + f_{x_t}' {V}_{t+1} f_{x_t} \nonumber\\
- & (f_{x_t}' {V}_{t+1} f_{u_t} + \bar{\lambda}_{t+1}f_{x_tu_t})' (R+f_{u_t}' {V}_{t+1} f_{u_t})^{-1} \nonumber \\
&(f_{u_t}' {V}_{t+1} f_{x_t} + \bar{\lambda}_{t+1}f_{x_tu_t}), \\
k_t &= (R+f_{u_t}' {V}_{t+1} f_{u_t})^{-1} (R\bar{u}_t + f_{u_t}' v_{t+1}), \\
K_t &= (R+f_{u_t}' {V}_{t+1} f_{u_t})^{-1} (f_{u_t}' {V}_{t+1} f_{x_t} + \bar{\lambda}_{t+1}f_{x_tu_t}),
\end{align*}
while the solution to the ILQR \eqref{N-LQR} is given by::
\begin{align*}
v_t &= \ell_{x_t} + f_{x_t}' v_{t+1} - f_{x_t}' V_{t+1} f_{u_t} (R+f_{u_t}' V_{t+1} f_{u_t})^{-1} (f_{u_t}' v_{t+1} + {R}\bar{u}_t), \\
V_t &= \ell_{x_t,x_t} + f_{x_t}' V_{t+1} f_{x_t} - f_{x_t}' V_{t+1} f_{u_t} (R+f_{u_t}' V_{t+1} f_{u_t})^{-1} f_{u_t}' V_{t+1} f_{x_t},  \\
k_t &= (R+f_{u_t}' V_{t+1} f_{u_t})^{-1} (R\bar{u}_t + f_{u_t}' v_{t+1}), \\
K_t &= (R+f_{u_t}' V_{t+1} f_{u_t})^{-1} (f_{u_t}' V_{t+1}f_{x_t}).
\end{align*}

Then, using the fact that iLQR is a special case of the modified SQP \eqref{M-QP}, it follows that iLQR always results in a descent direction and a suitable line search assures us that the next iterate is also feasible. Further, using the analogy of the Newton LQR \eqref{N-LQR} with the Newton SQP \eqref{N-QP}, it follows that the Newton LQR need not result in a descent direction, especially when far from an optimum. \\
Finally, a note regarding the line search in ILQR is due here. Typically, the updated control $\delta u_t$ is applied to the original nonlinear system along the nominal to generate a new nominal trajectory: $x_{t+1} = \bar{x}_t + \delta x_{t+1} = f(\bar{x}_t + \alpha \delta x_{t}, \bar{u}_t + \delta u_t)$ starting at $\alpha = 1$ and with $\delta x_0 = 0$. The solution is accepted if the new cost satisfies the following condition:
\begin{align*}
    \frac{J_{k+1} - J_k}{\alpha d_x' \nabla J_k} > \sigma > 0,
\end{align*}
for some user-defined $\sigma$, where $J_k$, $\nabla J_k$, represent the cost/ cost gradient at iteration $k$, while $J_{k+1}$ represents the cost under the updated control sequence. Note that the denominator above represents the cost change if linearity holds. If the inequality above is not satisfied, $\alpha$ is reduced till the inequality is satisfied (via a line search): that this is always possible is due to the fact that the direction $d_x$ is a cost descent direction. In essence, the feasibility is maintained by generating the new trajectory using the nonlinear dynamics, while the line-search assures us that we reach a better solution in terms of the cost. The development above can be summarized as follows.\\

\begin{proposition}
    The solution to the ILQR problem is always a descent direction for the cost $c(x,u)$ and the next solution can be made feasible via a suitable line-search procedure. In contrast, it is not necessary that the solution to the Newton LQR be a descent direction.\\
\end{proposition}

\begin{remark}
    Albeit paucity of space does not allow us to show the convergence proof for iLQR, it follows under standard assumptions (see \cite{Wang2025SearchFeedbackRL}).
\end{remark}

\begin{remark}
   In reference \cite{Pantoja_SN}, the author proposes a stagewise Newton procedure to find the Newton direction for the optimal control problem posed in this paper. In that work, the author takes a terminal cost/ Mayer approach to the problem and poses the problem as an unconstrained problem, unlike the constrained SQP route taken in this work. A careful look at the equations in that paper clearly shows that they are not the same as the equations we have for Newton LQR above. In particular, rather than the variables $\bar{\lambda}_t$ being the multipliers from the previous iteration, they are replaced by the adjoint variables $\nu_t$ from the current iteration which follow the equation: ${\nu}_t = l_{x_t} + f_{x_t}'{\nu}_{t+1}$, i.e., the calculation is equivalent to the first equation in the Newton LQR solution, except we do not have the $-(f_{x_t}' {V}_{t+1} f_{u_t} + \bar{\lambda}_{t+1}f_{x_tu_t})' (R+f_{u_t}' {V}_{t+1} f_{u_t})^{-1} (f_{u_t}' v_{t+1} + {R}\bar{u}_t)$ term after the first two terms. Noting that the term $(R\bar{u}_t + f_{u_t}' v_{t+1}) \approx 0$ close to an optimum, this clearly suggests that the approximation in \cite{Pantoja_SN} holds only close to an optimum, whereas the Newton LQR equations hold far from the optimum.
\end{remark}
\subsection{The Case of DDP}
In DDP, we start with the Bellman equation:
\begin{equation}
    V_t^*(x_t) = \min_{u_t} \left[ c(x_t, u_t) + V_{t+1}^*\left( f(x_t, u_t) \right) \right] \tag{B}
\end{equation}

and it is assumed that given a nominal $(\bar{x}_t, \bar{u}_t)$, the above equation can be written in the perturbation form:
\begin{equation}
    \Delta V_t^*(x_t) = \min_{\delta u_t} \left[ \Delta c(\delta x_t, \delta u_t) + \Delta V_{t+1}^*(\delta x_{t+1}) \right] \tag{PB}
\end{equation}
where $(\delta x_t, \delta u_t)$ are the perturbations from $(\bar{x}_t, \bar{u}_t)$, which is then expanded to second order to obtain an improved nominal $(\bar{x}_t', \bar{u}_t')$ \cite{JacobsonMayne1970DDP}.

However, (PB) does not follow from (B) unless $\bar{u}_t = u^*(\bar{x}_t)$, the optimal control action at nominal state $\bar{x}_t$. (B) may be written as:
\begin{equation}
    V_t^*(x_t) = \min_{u_t} Q(x_t, u_t)
\end{equation}
where
\begin{equation}
    Q(x_t, u_t) = c(x_t, u_t) + V_{t+1}^* \left( f(x_t, u_t) \right).
\end{equation}

Let $u_t^*(\cdot)$ represent the optimal feedback law at time $t$.
Then:
\begin{equation}
    V_t^*(x_t) = Q(x_t, u_t^*(x_t)).
\end{equation}
and doing a perturbation expansion around $\bar{x}_t$, one gets:
\begin{equation}
    V_t^*(\bar{x}_t + \delta x_t) = Q \left( \bar{x}_t + \delta x_t, \; u_t^*\left( \bar{x}_t + \delta x_t \right) \right).
\end{equation}
Therefore,
\begin{align*}
 &V_t^*(\bar{x}_t + \delta x_t) = Q(\bar{x}_t, \bar{u}_t^*) + \left. \frac{\partial Q}{\partial x_t} \right|_{\bar{x}_t} \delta x_t + \left. \frac{\partial Q}{\partial u_t} \right|_{u_t^*} \delta u_t + \text{H.O.T} \\
&= Q_t^* + \left(\left. \frac{\partial Q}{\partial x_t} \right|_{\bar{x}_t} + \left. \frac{\partial Q}{\partial u_t} \right|_{u_t^*} \cdot \left.\frac{d u_t^*}{d x_t} \right|_{\bar{x}_t}\right) \delta x_t + \text{H.O.T} \tag{\#} 
\end{align*}
Now, consider some nominal $(\bar{x}_t, \bar{u}_t)$ where $\bar{u}_t \neq u_t^*$.
Then, an expansion of $\overline{V}_t(x_t)$ leads to:
\begin{align*}
\overline{V}_t (\bar{x}_t + \delta x_t) = \min_{\delta u_t} \Big[ Q(\bar{x}_t, \bar{u}_t) + \Delta Q(\delta x_t, \delta u_t) \Big]
\end{align*}
However, similar to the optimal case, in general,
\begin{align*}
\bar{u}_t(x_t) = \overline{u}_t(\bar{x}_t) + \frac{d \overline{u}_t}{dx_t}\Big|_{\bar{x}_t} \delta x_t + \text{H.O.T.}
\end{align*}
where $\overline{u}_t(\bar{x}_t) = \bar{u}_t$.
Then:
\begin{align}
\overline{V}_t (\bar{x}_t + \delta x_t) = \bar{Q}_t + \left( \left. \frac{\partial Q}{\partial x_t} \right|_{\bar{x}_t}  + \left.\frac{\partial Q}{\partial u_t} \right|_{\bar{u}_t} \frac{d \overline{u}_t}{d x_t}\Big|_{\bar{x}_t}  \right)\delta{x}_t + \text{H.O.T.} \tag{\#\#}
\end{align}
Note that $\bar{Q}_t \ne Q_t^*$, $u_t^* \neq \overline{u}_t$, $\left. \frac{d u_t^*}{d x_t} \right|_{\bar{x}_t} \neq \frac{d \overline{u}_t}{d x_t} \Big|_{\bar{x}_t}$,
and $\frac{\partial Q}{\partial u_t}\Big|_{\bar{u}_t} \ne \frac{\partial Q}{\partial u_t}\Big|_{u_t^*}$.
Therefore, from (\#) and (\#\#), it follows that:
\[
\overline{V}_t(\bar{x}_t + \delta x_t) \neq V_t^*(\bar{x}_t + \delta x_t),
\]
and the perturbation equation (PB) cannot follow from the DP equation (B) unless $\bar{u}_t = u_t^*(\bar{x}_t)$. A final note here is that the optimal feedback law $u_t^*(.)$ has to be smooth enough for the perturbation expansions to hold. \\
The above development can be summarized as the following result.
\begin{proposition}
    A perturbation expansion of the DP equation is feasible if and only if the perturbation is with respect to an optimal nominal path $(x_t^*,u_t^*)$ and given the feedback law $u_t^*(.)$ is smooth enough.\\
\end{proposition}

Thus, this implies that the DDP expression arising from the second order approximation of (PB) is essentially a heuristic approximation of the Dynamic Programming equation (B).
The DDP equations can be seen to be identical to the Newton LQR equation except that in DDP, the multipliers $\bar{\lambda}_t$ from the previous iteration are replaced by the current nominal multiplier component $v_t$. In particular, when close to the optimum, $v_t \approx \bar{\lambda}_t$, and thus, DDP can be construed as a heuristic approximation to the Newton LQR equations, at least close to an optimum. Also, note that, owing to this interpretation, DDP and the stagewise Newton approach in \cite{Pantoja_SN}, are essentially identical close to an optimum. Thus, it follows that the DDP solution cannot be assured to be a descent direction for the cost function when far from an optimum.




\section{Contrasting the algorithmic performance of DDP and iLQR}

While the preceding discussion outlined the theoretical aspects of iLQR and DDP, it is equally important to understand how these algorithms behave in practice. In this section, we present a direct comparison of their performance on representative control tasks, explicitly focusing on the unregularized forms of both iLQR and DDP. By analyzing convergence rates, predicted versus actual cost reductions, and step sizes, we aim to reveal the practical advantages and shortcomings that arise solely from their respective local modeling choices, independent of any stabilization or damping heuristics.

\subsection{Cost-Descent Direction}
In this subsection, we examine the descent directions predicted by iLQR and DDP. Under mild assumptions, iLQR consistently produces a descent direction, whereas DDP does not offer the same guarantee.\\
The net ``expected cost reduction" in these algorithms is given by:
\begin{equation}
    \Delta J = -\Bigg{(} \alpha - \frac{\alpha^2}{2} \Bigg{)}\sum_{t=0}^{T-1} k_{t}^{T}\underbrace{(R+f_{u_t}' {V}_{t+1} f_{u_t})^{-1}}_{Q_{uu}^{-1}}k_{t},
\end{equation}
where all variables follow the same definitions as in the preceding discussion, and in particular, $\alpha$ is the line search parameter. It should be noted that the term $Q_{uu}$ should be positive definite for the reduction in cost. Under some mild assumptions, this term is always positive definite for iLQR. However, it may become negative-definite/indefinite for DDP. We consider the cartpole swing up and pendulum swing up task to show this comparison.

\begin{figure}[!htbp]
    \centering
    \subfloat[Cart-pole swing-up.\label{subfig:cartpole}]{\includegraphics[width=0.8\linewidth]{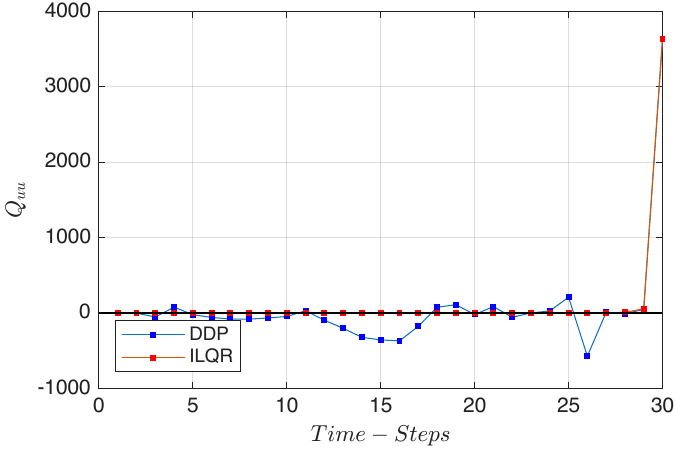}}
   \hfill
   \subfloat[Pendulum swing-up.\label{subfig:pendulum}] {\includegraphics[width=0.8\linewidth]{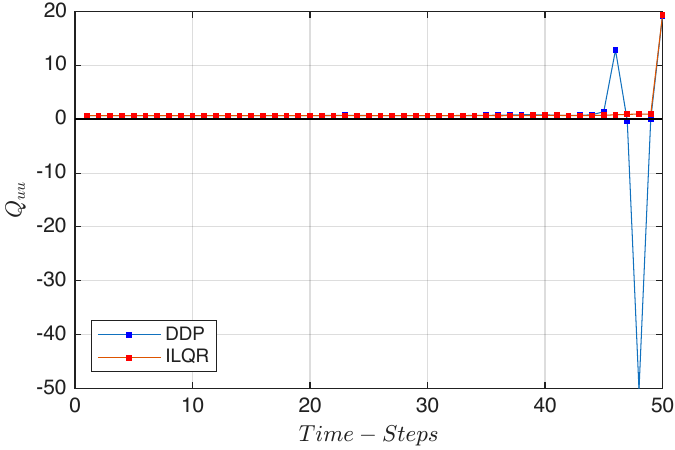}}%
    \caption{The plot of $Q_{uu}$ for the cartpole and pendulum problems for the first iteration starting with random initial guesses.}
    \label{fig:Quu_pd}
\end{figure}

Figure \ref{fig:Quu_pd} illustrates that $Q_{uu}$ remains positive in the cartpole and pendulum problems for iLQR. Since there is only a single control input, $Q_{uu}$ is scalar in these cases. However, in DDP, it becomes negative at certain points, which can lead to an increase in cost at some time steps. This mixture of positive and negative values for $Q_{uu}$ across the trajectory necessitates strong regularization, as the predicted change in cost may otherwise increase. For the pendulum case shown in Fig.~\ref{subfig:pendulum}, the predicted cost variations are $\Delta J_{DDP} = -224.0154$ and $\Delta J_{iLQR} = -272.9982$. Both methods predict a cost decrease, with iLQR performing better in this unregularized setting since $Q_{uu}$ remains positive throughout. In contrast, Fig.~\ref{subfig:cartpole} presents the $Q_{uu}$ profile for the cart-pole. Here, DDP yields $\Delta J_{DDP} = 3.0172 \times 10^7$, indicating a cost increase, while iLQR predicts $\Delta J_{iLQR} = -5.5692 \times 10^5$, reflecting a cost reduction.

\subsection{Cooling of the line-search parameter $\alpha$}
It is a common belief that DDP accepts larger steps than iLQR; equivalently, for a given decrease in cost, the line-search parameter~$\alpha$ is often higher under DDP than under iLQR, motivated by the idea that DDP forms a better local model by including a second-order expansion of the dynamics. In this section, we show empirical evidence that DDP can lead to excessive cooling of the line-search parameter, $\alpha$, and hence, can take a much smaller step than iLQR.

\begin{figure}[!htbp]
    \centering
    \subfloat[Learning rate vs no. of iterations.\label{subfig:alpha_pendulum}]{\includegraphics[width=.48\linewidth]{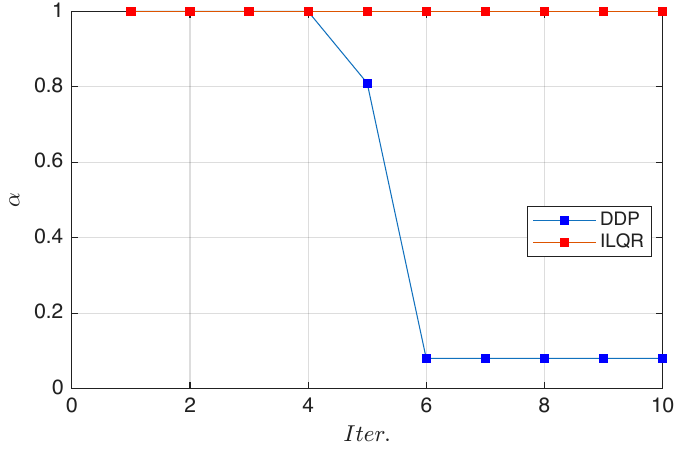}}
   \hfill
   \subfloat[Cost vs no. of iterations.\label{subfig:cost_pendulum}] {\includegraphics[width=0.48\linewidth]{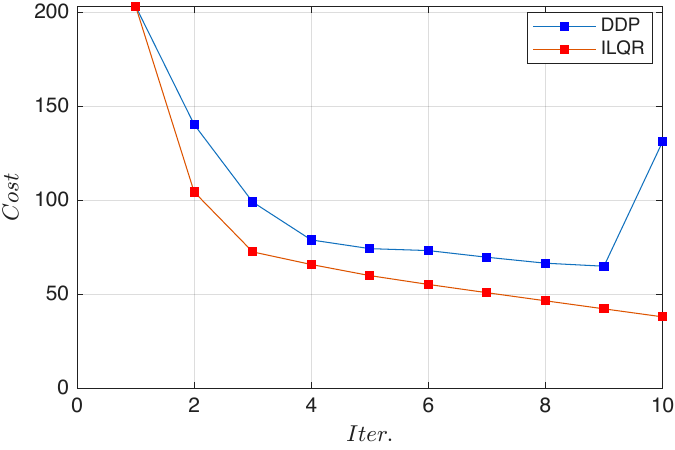}}%
    \caption{The plot of learning rate $\alpha$ and cost vs iteration for the pendulum swing-up task starting with random initial guesses. DDP cools down $\alpha$ as the cost predictions are corrupted.}
    \label{fig:cool_pendulum}
\end{figure}

\begin{figure}[!htbp]
    \centering
    \subfloat[Learning rate vs no. of iterations.\label{subfig:alpha_cartpole}]{\includegraphics[width=.48\linewidth]{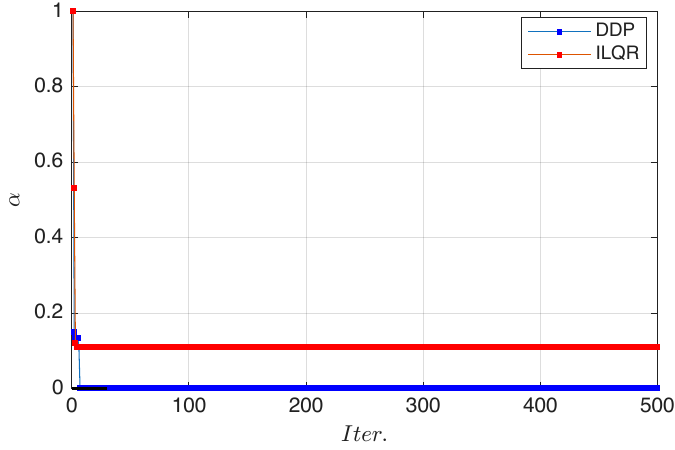}}
   \hfill
   \subfloat[Cost vs no. of iterations.\label{subfig:cost_cartpole}] {\includegraphics[width=0.48\linewidth]{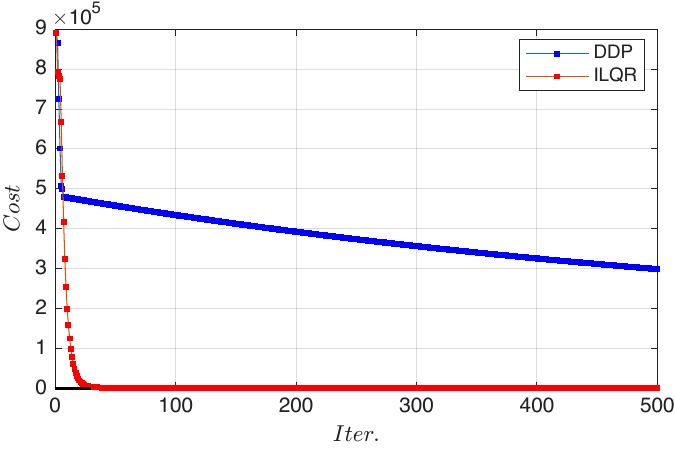}}%
    \caption{The plot of learning rate $\alpha$ and cost vs iteration for the cart-pole problem starting with random initial guesses. DDP cools down $\alpha$ as the cost change predictions are corrupted.}
    \label{fig:cool_cartpole}
\end{figure}

iLQR ties both prediction and realization to a single (first-order) state-propagation model, which bounds the necessary backtracking. DDP, by mixing linear and second-order information introduces additional $\alpha$-sensitivity that can cause excessive cooling of the line-search parameter.

Figures \ref{fig:cool_pendulum} and \ref{fig:cool_cartpole} contrast the variation of cost and learning rate, $\alpha$, for iLQR and DDP over successive iterations. As seen in Fig.~\ref{fig:cool_pendulum},  for the pendulum swing up, the learning rate remains fixed at $1$ for iLQR, while for DDP it drops to smaller values after iteration 4. Even during the first four iterations, when $\alpha = 1$ for both methods, iLQR achieves lower costs than DDP. In addition, DDP exhibits a noticeable cost increase from iteration 9 to 10, demanding a need for regularization. Likewise, for cartpole swing up, in Fig.~\ref{fig:cool_cartpole}, the learning rate for DDP decreases to nearly zero after a few iterations, resulting in a very slow cost reduction in subsequent steps. By contrast, iLQR maintains a learning rate of about $0.1$ and reduces the cost more effectively.

\begin{figure}[!htbp]
    \centering
    \subfloat[Pendulum swing-up.\label{subfig:comp_pendulum}]{\includegraphics[width=.80\linewidth]{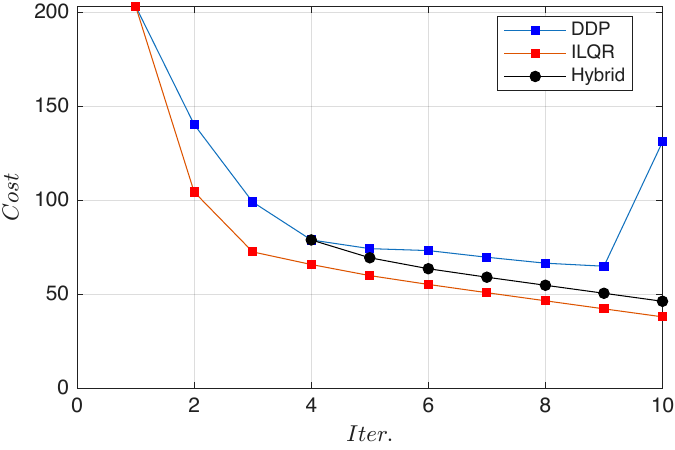}}
   \hfill
   \subfloat[Cart-pole swing-up.\label{subfig:comp_cartpole}] {\includegraphics[width=0.80\linewidth]{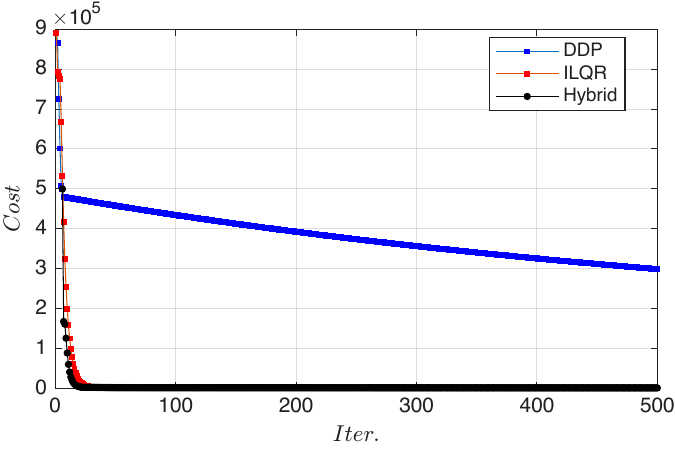}}%
    \caption{A comparison of iLQR, DDP and hybrid solution. The hybrid solution depicts iLQR applied to the point where DDP solution slows down as $\alpha$ reduces drastically (depicted in Figures \ref{fig:cool_pendulum} and \ref{fig:cool_cartpole}).}
    \label{fig:correction_lr}
\end{figure}

We applied iLQR at the point where DDP significantly reduces the line-search parameter. The resulting performance is shown in Fig.~\ref{fig:correction_lr}, where it is labeled as ``Hybrid.'' The plots clearly indicate that switching to iLQR at the stage where DDP cools down improves performance. Unlike DDP, iLQR maintains a stable learning rate and shows better convergence for both the cart-pole and pendulum swing-up tasks.

\subsection{DDP may predict an unrealistic cost decrease}
This subsection highlights that the DDP quadratic model can predict a decrease in cost that is not physically attainable, whereas the iLQR model does not. Formally,
\begin{align}
    \overline{J} + \Delta J_{\mathrm{iLQR}} \;\ge\; J_{\min}, \qquad
    \overline{J} + \Delta J_{\mathrm{DDP}} \;<\; J_{\min},
\end{align}
where $\overline{J}$ is the cost using the current policy (nominal control input). $J_{min}$ is the minimum value by which the cost is bounded. In our case, we have considered a purely quadratic cost with no cross-coupling between states and control. So, our cost is bounded below by zero for each problem. 

\begin{table}
    \centering
    \begin{tabular}{|c|c|c|c|}
    \hline
      Iter.   & $J$ & $\Delta J$ & $J_{\text{pred}} = J + \Delta J$\\
      \hline
        0 & $701.4661$ & $-377.7732$ & $323.6929$ \\
        \hline
        2 & $149.8840$ & $-51.4378$ & $98.4462$ \\
        \hline
        5 & $2.872408 \times 10^{5}$ & $-3.3685 \times 10^{5}$ & $-4.9609 \times 10^{4}$ \\
    \hline
    \end{tabular}
    \caption{Change in cost predicted by DDP over various iterations starting from random initial guesses for the pendulum swing-up task.}
    \label{tab:cost_pred_DDP_pen}
\end{table}

\begin{table}
    \centering
    \begin{tabular}{|c|c|c|c|}
    \hline
      Iter.   & $J$ & $\Delta J$ & $J_{\text{pred}} = J + \Delta J$\\
      \hline
        0 & $8.897408 \times 10^{5}$ & $-8.1926 \times 10^{5}$ & $7.0481 \times 10^{4}$ \\
        \hline
        3 & $6.011446 \times 10^{5}$ & $-9.3269e \times 10^{4}$ & $5.0788 \times 10^{5}$ \\
        \hline
        5 & $4.998068 \times 10^{5}$ & $-3.6467 \times 10^{7}$ & $-3.5968 \times 10^{7}$ \\
        \hline
    \end{tabular}
    \caption{Change in cost predicted by DDP over various iterations starting from random initial guesses for the cart-pole swing-up task.}
    \label{tab:cost_pred_DDP}
\end{table}

Table \ref{tab:cost_pred_DDP} indicates that at iteration 5, the cost predicted by DDP becomes negative ($-3.6\times 10^{7}$). Since the experiment used a purely quadratic cost with a lower bound of $0$ ($J_{\text{min}} = 0$), such a negative prediction is not feasible. This issue does not occur with iLQR.

A similar problem is shown in Table \ref{tab:cost_pred_DDP_pen}, which shows the same problem for the pendulum swing-up task. DDP predicts a change in cost $\Delta J$ such that $J + \Delta J < J_{min}$. It should also be noted that the cost at iteration 5 is more than that of previous iterations, as we are considering unregularized DDP, and it tends to increase the cost because of $Q_{uu}$ being non-positive definite. 

\subsection{Faster convergence of DDP over iLQR near optima}
In certain scenarios, DDP demonstrates faster convergence than iLQR. In particular, near the optimal solution, a DDP step effectively corresponds to a Newton step. As a result, DDP often takes larger and more effective steps in proximity to the optimum, leading to quicker convergence in many cases.

\begin{figure}[!htbp]
    \centering
    \subfloat[Pendulum swing-up.\label{subfig:DDP_pen}]{
    \includegraphics[width=.45\linewidth]{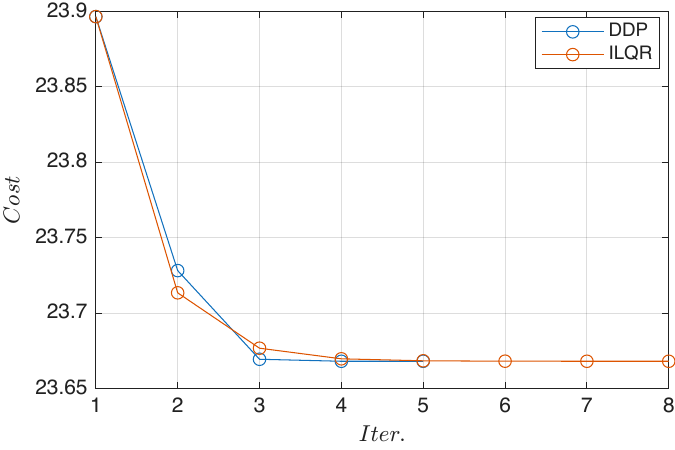}
    }
    \hfill
    \subfloat[Cart-pole swing-up\label{subfig:DDP_cart}]{
    \includegraphics[width=0.45\linewidth]{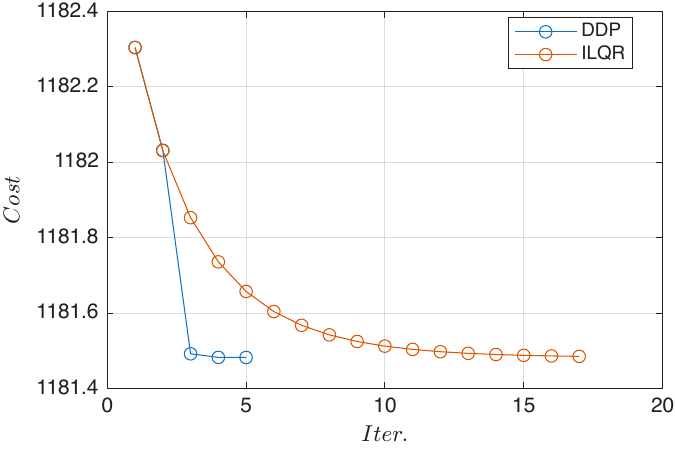}
    }
    \caption{Cost vs. number of iterations. With an initial control guess near the optimum, DDP reaches the optimum in fewer iterations.}
    \label{fig:DDP_fast}
\end{figure}

Figure \ref{fig:DDP_fast} illustrates that when initialized sufficiently close to the optimum, DDP tends to converge in fewer iterations by taking larger steps. In Fig.~\ref{subfig:DDP_pen}, for the pendulum swing-up task, DDP converges within 5 iterations, whereas iLQR requires 8. Similarly, in Fig.~\ref{subfig:DDP_cart}, for the cartpole swing-up, DDP converges in 5 iterations compared to 17 iterations for iLQR. Thus, in the neighborhood of an optimal solution, DDP typically achieves faster convergence than iLQR.



\section{Conclusion}

We have studied the solution to optimal control problems, in particular, Newton's approach, Differential Dynamic Programming (DDP), and iterative Linear Quadratic Regulation (iLQR), through a unified constrained SQP framework. The critical insight is that the Newton step in a constrained optimization is equivalent to a QP, and consequently, an LQR problem when solving an optimal control problem. The analysis reveals that, although DDP can be construed as an approximation to the Newton step close to an optimum, it often leads to inconsistent cost predictions and unreliable updates, especially when far from an optimum. In contrast, iLQR provides consistent descent directions and admits theoretical convergence guarantees, making it far more robust in practice. Experiments on the pendulum and cartpole tasks confirm that iLQR delivers stable performance, whereas DDP may produce dubious solutions unless heavily regularized, which means there are no guarantees of convergence when far from an optimum. These results also establish iLQR as a principled SQP approach to optimal control rather than a mere approximation of DDP by neglecting the second-order terms, and highlight the importance of consistency over second-order information in optimal control problems. In fact, our theoretical and empirical results show that not including the second-order information is actually critical to the performance, thereby showing that doing less is actually more in this instance.

\printbibliography

@book{JacobsonMayne1970DDP,
  author    = {David H. Jacobson and David Q. Mayne},
  title     = {Differential Dynamic Programming},
  series    = {Modern Analytic and Computational Methods in Science and Mathematics},
  volume    = {24},
  address   = {New York},
  publisher = {American Elsevier},
  year      = {1970},
  pages     = {208},
  isbn      = {0444000704}
}

@inproceedings{todorov2005ilqg,
  title={A generalized iterative LQG method for locally-optimal feedback control of nonlinear stochastic systems},
  author={Todorov, Emanuel and Li, Weiwei},
  booktitle={Proceedings of the American Control Conference},
  pages={300--306},
  year={2005}
}

@article{tassa2012synthesis,
  title={Synthesis and stabilization of complex behaviors through online trajectory optimization},
  author={Tassa, Yuval and Erez, Tom and Todorov, Emanuel},
  journal={Proceedings of the IEEE/RSJ International Conference on Intelligent Robots and Systems},
  pages={4906--4913},
  year={2012}
}

@article{tassa2014control,
  title={Control-limited differential dynamic programming},
  author={Tassa, Yuval and Mansard, Nicolas and Todorov, Emanuel},
  journal={Proceedings of the IEEE International Conference on Robotics and Automation},
  pages={1168--1175},
  year={2014}
}

@article{levine2013guided,
  title={Guided policy search},
  author={Levine, Sergey and Koltun, Vladlen},
  journal={Proceedings of the 30th International Conference on Machine Learning},
  pages={1--9},
  year={2013}
}

@article{manchester2017dircol,
  title={Dircol: Direct collocation methods for trajectory optimization},
  author={Manchester, Zachary and Kuindersma, Scott},
  journal={IEEE Robotics and Automation Magazine},
  volume={24},
  number={3},
  pages={36--45},
  year={2017}
}

@article{mayne1965ddp,
  title={A second-order gradient method for determining optimal trajectories of non-linear discrete-time systems},
  author={Mayne, David Q},
  journal={International Journal of Control},
  volume={3},
  number={1},
  pages={85--95},
  year={1965}
}

@article{boggs1995sqp,
  title={Sequential quadratic programming methods for large-scale nonlinear optimization},
  author={Boggs, Paul T. and Tolle, Jon W.},
  journal={Mathematical Programming},
  volume={62},
  number={1-3},
  pages={3--30},
  year={1995},
  publisher={Springer}
}

@article{Pantoja_SN,
author = {J. F. A. De O. Pantoja},
title = {Differential dynamic programming and Newton's method},
journal = {International Journal of Control},
volume = {47},
number = {5},
pages = {1539--1553},
year = {1988},
publisher = {Taylor \& Francis},
doi = {10.1080/00207178808906114}}

@article{NgangaWensing2021RAL,
  author  = {John N. Nganga and Patrick M. Wensing},
  title   = {Accelerating Second-Order Differential Dynamic Programming for Rigid-Body Systems},
  journal = {IEEE Robotics and Automation Letters},
  year    = {2021},
  volume  = {6},
  number  = {4},
  pages   = {7659--7666},
  doi     = {10.1109/LRA.2021.3098928}
}

@inproceedings{Plancher2017CUDP,
  author    = {Brian Plancher and Zachary Manchester and Scott Kuindersma},
  title     = {Constrained Unscented Dynamic Programming},
  booktitle = {2017 IEEE/RSJ International Conference on Intelligent Robots and Systems (IROS)},
  year      = {2017},
  month     = sep,
  pages     = {5674--5680},
  address   = {Vancouver, BC, Canada},
  publisher = {IEEE},
  doi       = {10.1109/IROS.2017.8206457}
}

@techreport{LiaoShoemaker1992,
  author      = {Li-Zhi Liao and Christine A. Shoemaker},
  title       = {Advantages of Differential Dynamic Programming over Newton's Method for Discrete-Time Optimal Control Problems},
  institution = {Cornell University},
  type        = {Technical Report},
  number      = {ctc92tr97},
  address     = {Ithaca, NY},
  year        = {1992},
  url         = {https://ecommons.cornell.edu/items/2e4fec06-1f95-4576-b269-7660105605c8}
}

@article{Wang2025SearchFeedbackRL,
  author  = {Ran Wang and Aayushman Sharma and Karthikeya S. Parunandi and Raman Goyal and Mohamed Naveed Gul Mohamed and Suman Chakravorty},
  title   = {The Search for Feedback in Reinforcement Learning},
  journal = {Journal of Dynamic Systems, Measurement, and Control},
  year    = {2025},
  volume  = {147},
  number  = {6},
  pages   = {061002},
  doi     = {10.1115/1.4068705}
}

@misc{zhang2025wholebodymodelpredictivecontrollegged,
      title={Whole-Body Model-Predictive Control of Legged Robots with MuJoCo}, 
      author={John Z. Zhang and Taylor A. Howell and Zeji Yi and Chaoyi Pan and Guanya Shi and Guannan Qu and Tom Erez and Yuval Tassa and Zachary Manchester},
      year={2025},
      eprint={2503.04613},
      archivePrefix={arXiv},
      primaryClass={cs.RO},
      url={https://arxiv.org/abs/2503.04613}, 
}

@article{Aziz2019HDDP_CR3BP,
  author  = {Jonathan D. Aziz and Daniel J. Scheeres and Gregory Lantoine},
  title   = {Hybrid Differential Dynamic Programming in the Circular Restricted Three-Body Problem},
  journal = {Journal of Guidance, Control, and Dynamics},
  year    = {2019},
  volume  = {42},
  number  = {5},
  pages   = {963--975},
  doi     = {10.2514/1.G003617},
  url     = {https://doi.org/10.2514/1.G003617}
}

@article{Russell1,
author= {Lantoine,Gregory and Russell,Ryan P.},
year={2012},
month={08},
title={A Hybrid Differential Dynamic Programming Algorithm for Constrained Optimal Control Problems. Part 1: Theory},
journal={Journal of Optimization Theory and Applications},
volume={154},
number={2},
pages={382-417},
language={English},
}

@article{Russell2,
author={Lantoine,Gregory and Russell,Ryan P.},
year={2012},
month={08},
title={A Hybrid Differential Dynamic Programming Algorithm for Constrained Optimal Control Problems. Part 2: Application},
journal={Journal of Optimization Theory and Applications},
volume={154},
number={2},
pages={418-442},
isbn={00223239},
language={English},
}

@inproceedings{Oshin2022PDDP_RSS,
  author    = {Alex Oshin and Matthew D. Houghton and Michael J. Acheson and Irene M. Gregory and Evangelos A. Theodorou},
  title     = {Parameterized Differential Dynamic Programming},
  booktitle = {Robotics: Science and Systems XVIII (RSS)},
  year      = {2022},
  address   = {New York City, NY, USA},
  month     = jun,
  url       = {https://www.roboticsproceedings.org/rss18/p046.pdf}
}

@techreport{Noyes2021RobustEntryGuidance,
  author      = {Connor David Noyes},
  title       = {Robust Optimal Entry Guidance for Future Mars Landers},
  year        = {2021},
  institution = {University of California, Irvine},
  type        = {Technical Report},
  url         = {https://escholarship.org/uc/item/7bv3q5ws},
  note        = {Open-access UC eScholarship report}
}

\addtolength{\textheight}{-12cm}   








\end{document}